\documentclass[letterpaper, 10 pt, conference]{ieeeconf}  

\usepackage[T1]{fontenc}         
\usepackage{graphicx}
\usepackage{caption}
\usepackage{subcaption} 
\usepackage{float}   
\usepackage{fancybox}		  
\usepackage{listings} 
\usepackage{amsfonts}
\usepackage{amsmath}
\usepackage{flushend}

\newtheorem{theorem}{Theorem}[section]

\newcommand{\R}{\mathbb{R}}

\IEEEoverridecommandlockouts                              

\title{\LARGE \bf
Developmental Partial Differential Equations
}

\author{Nastassia Pouradier Duteil$^{1}$,  Francesco Rossi$^{2}$,
Ugo Boscain$^{3}$, Benedetto Piccoli$^{1}$ 
\thanks{$^{1}$ Department of Mathematical Sciences, Rutgers University - Camden, Camden, NJ.
        {\tt\small piccoli@camden.rutgers.edu, nastassia.pouradierduteil@rutgers.edu}}%
\thanks{$^{2}$ Aix Marseille Universit\'e, CNRS, ENSAM, Universit\'e de Toulon, LSIS UMR 7296, 13397, Marseille, France. 
        {\tt\small francesco.rossi@lsis.org}}
\thanks{$^{3}$ CNRS, CMAP, Ecole Polytechnique, Palaiseau, France, \& Team GECO, INRIA Saclay. {\tt\small ugo.boscain@polytechnique.edu}}%
}

\begin{document}

\maketitle
\thispagestyle{empty}
\pagestyle{empty}

\begin{abstract}
In this paper, we introduce the concept of {\em Developmental Partial Differential Equation} (DPDE), which consists of a Partial Differential Equation (PDE) on a time-varying manifold with complete coupling between the PDE and the manifold's evolution. In other words, the manifold's evolution depends on the solution to the PDE, and vice versa the differential operator of the PDE depends on the manifold's geometry. DPDE is used to study a diffusion equation with source on a growing surface whose growth depends on the intensity of the diffused quantity. The surface may, for instance, represent the membrane of an egg chamber and the diffused quantity a protein
activating a  signaling pathway leading to growth.
Our main objective is to show controllability of the surface shape using a fixed source with variable intensity for the diffusion. More specifically, we look for a control driving a symmetric manifold shape to any other symmetric shape in a given time interval. For the diffusion we take directly the Laplace-Beltrami operator of the surface, while the surface growth is assumed to be equal to the value of the diffused quantity. 
We introduce a theoretical framework, provide approximate controllability and show numerical results. Future applications include a specific model for the oogenesis of \textit{Drosophila melanogaster}.
\end{abstract}

\section*{INTRODUCTION}
Many mathematical models aim to reproduce biological mechanisms, including those involving growth of living organisms (see \cite{Chauvet, DeV, Murray, Segel}). 
This concerns the field of Developmental Biology, which is devoted to the study of growth and development of organisms, as well as the genetic control of cell growth, differentiation and morphogenesis (see \cite{Spradling, Wolpert}). 

We focus on providing a mathematical framework for understanding growth
of organisms induced by signaling pathways. 
Many approaches for modeling biological growth
were proposed in the literature, see for instance \cite{Jones} and references therein.
In that case, a tissue is either regarded as a continuum or as a collection of cells. The latter microscopic approach is based on discrete models such as cellular automata. On the other hand, Partial Differential Equations are natural
for the former macroscopic point of view, where the dependent variable 
usually represents mass concentration. In most models the growth of mass
is then assigned as an internal mechanism, e.g. proportional to the mass
concentration, or as an external one. The point of view is similar to that of elasticity,
but the equations are modified to include growth (so violating conservation
of mass).

There exist abundant investigations on Turing Patterns (TP), which are generated by two chemicals due to instabilities caused by diffusion (see \cite{Turing, Varea}).
Papers also investigated effects on stability, geometry and growth in TP due to growth rate, curvature and other characteristics of chemical interactions (\cite{Baker, Chrisholm, Crampin, Gjorgjieva, Lefevre, Miura, Plaza, Simpson}). However, to our knowledge, there is no clear experimental evidence of model organisms showing patterns via mechanisms predicted by the corresponding mathematical model.

In our approach, we aim at having a faithful representation of the fact that growth
is mostly regulated by signaling pathways.
In order to achieve this,
we design a mathematical framework based on two main ingredients.
The first is a smooth topological manifold ${\cal M}$
(representing for instance an egg chamber, a tissue,
or even a single cell) which evolves in time. The second is a quantity $s$
evolving on the manifold according to a specific PDE.
More precisely, we assume that the manifold's change (usually growth)
depends directly on the quantity modeled by the PDE and, in turn,
the PDE operator depends on the (changing in time) manifold's geometry.
In order to define the model in a mathematically sound way,
we consider a manifold
${\cal M}$ embedded in a ambient Euclidean space $\R^n$
and parameterized (possibly locally) by variables $y\in\R^m$ via a map $\psi(t,y)$.
Moreover, the evolution of ${\cal M}$ is given by a vector field $v(s,x)$, $x\in\R^n$, depending on $s$. 
To mimic the genetic mechanisms regulating growth, we introduce
a control term which sets the intensity of a source for the quantity $s$. 
Finally the coupled system reads:
\begin{equation}\label{eq:universe}
\frac{\partial \psi (t,y)}{\partial t}=v(s(t,y),y),\qquad
X_{{\cal M}} s =S(u),
\end{equation}
where $X_{{\cal M}}$ is a differential operator defined on the manifold ${\cal M}$,
$S$ a source term and $u$ the control.\\
PDEs on manifolds form a vast subject, with a wide and varied literature. Previous works focusing on the control of PDEs on manifolds include \cite{Coron, Khapalov}.
Other works have dealt with many different
PDEs on manifolds. For instance, \cite{Struwe} focused
on harmonic maps on Riemannian surfaces, while the Cauchy problem for wave maps is treated in \cite{Shatah}. A general approach for hypoelliptic 
Laplacian on unimodular Lie groups can be found in \cite{Agrachev}.
Regarding controlled evolution on manifolds,
the Klein-Gordon equation on a 3-D compact manifold is considered in \cite{Laurent}.
All these works focus on PDEs on manifolds (and possibly their control) 
but, to our knowledge, our contribution is the first that deals with a completely
coupled system manifold-PDE of the type (\ref{eq:universe}).
Because of its biological meaning, we call such system a Developmental 
Partial Differential Equation (DPDE). 

The specific application we have in mind is that of understanding how the Gurken protein contributes to the formation of various \textit{Drosophila} eggshell structures. During oogenesis, Gurken is secreted near the oocyte nucleus, and then diffuses and is integrated by the epidermal growth factor receptor (EGFR), triggering a signaling pathway. 
Therefore, modeling evolution of the Gurken concentration during \textit{Drosophila} oogenesis is more complex than solving a PDE on a manifold due to the fact that the egg chamber grows throughout the process. 
Leveraging on natural symmetries of the egg chamber
we choose as ${\cal M}$ a one-dimensional symmetric manifold embedded
in $\R^2$ and initially equal to $S^1$. For the operator $X_{{\cal M}}$,
we simply choose the Laplace-Beltrami operator on ${\cal M}$. Finally,
the resulting DPDE is given by (\ref{eq:devpde}).\\
Our main aim is to show controllability in terms of the possible shapes
reachable from $S^1$ regulating one or more sources. We show
how to adapt the approach of Laroche, Martin and Rouchon \cite{Laroche},
proving flatness of the heat equation, to our setting and then
provide numerical studies.


The paper is organized as follows.
We begin by introducing the biological context that motivates our work, that is diffusion of Gurken in the growing \textit{Drosophila} egg chamber during oogenesis. 
We then establish the mathematical framework, starting with a general setting involving multiple reaction-diffusion equations on a dynamically evolving manifold. We simplify the problem by considering the diffusion of just one quantity (the growth signal) and its direct impact on the manifold's evolution. Finally, we show numerical simulations depicting various possible shapes that can be obtained by controlling one or two sources for the signal.



\section{BIOLOGICAL CONTEXT}

\textit{Drosophila melanogaster} is a commonly used model organism to study
cell signaling, tissue patterning and morphogenesis.
Tissue patterning and cell fate determination are guided by a handful of cell signaling pathways. 
For instance, the epidermal growth factor receptor (EGFR) signaling pathway controls many cell processes, including apoptosis (cell death) and cell migration. 
In particular for \textit{Drosophila} oogenesis, the TGF-alpha like ligand Gurken is secreted from near the oocyte nucleus,  
diffuses in the perivitelline space surrounding the oocyte, 
and signals through EGFR in the overlaying follicle cells. 
This will eventually gives rise to cell differentiation, forming various eggshell structures. 

These processes occur over the course of about 27 hours, during which the \textit{Drosophila} egg chamber undergoes morphological changes, including the follicle cells'
gradual movement over the oocyte (see Figure 
1.).
Thus, different cells are dynamically exposed to Gurken over time. Mathematically, this is equivalent to considering the oocyte nucleus as a moving source of Gurken. Moreover, during oogenesis the egg chamber grows and becomes wider and more elongated.

The diffusion of Gurken has been previously modeled to determine its impact on local genes (see \cite{Goentoro, Zartman}). These existing models consider diffusion at steady-state from a fixed nucleus, thus not accounting for the evolution of the egg chamber's shape.
 The aim of this work is to propose a mathematical setting that will allow to model the growth of the egg chamber as oogenesis progresses.
 
 \begin{figure}[h!]\label{fig:Zartman}
        \begin{center}
        \begin{subfigure}[b]{0.5\textwidth}
                \includegraphics[trim=0cm 13cm 0cm 0cm, clip=true, scale=0.35]{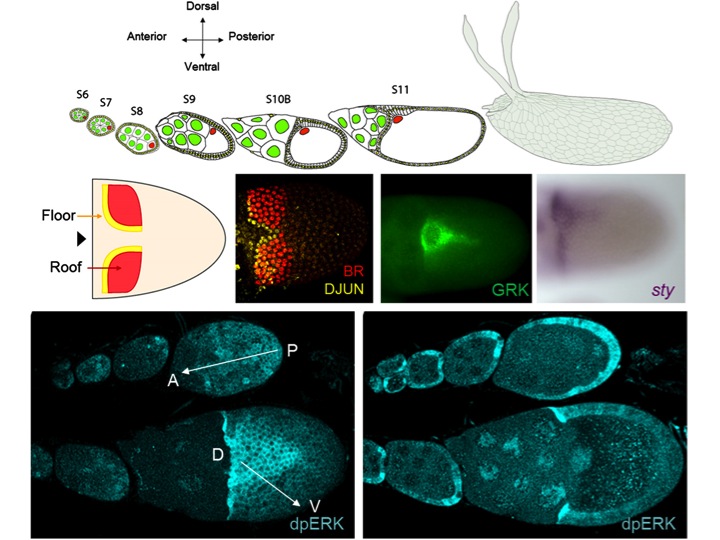}  
                \caption{\footnotesize Schematic of the eggshell growth during oogenesis from Stage 6 to Stage 11. The resulting egg with protruding dorsal appendages is shown on the right.}
        \end{subfigure}
                        \vspace{0.5cm}
        \begin{subfigure}[b]{0.5 \textwidth}
                \includegraphics[trim=0cm 0cm 0cm 11cm, clip=true, scale=0.35]{Zartman.jpg} 
                \caption{\footnotesize 
               The Gurken ligand diffuses in the oocyte and binds to EGFR, triggering a signaling cascade. The level of EGFR signaling can be assayed with dpERK antibody staining. The non-homogeneous activation of EGFR is responsible for the formation of various morphological features on Drosophila species' eggshells, such as the dorsal appendages or the dorsal ridge.
               }
        \end{subfigure}
        
        \caption{From \cite{Zartman} J.J. Zartman, L.S. Cheung, M. Niepielko, C. Bonini, B. Haley, N. Yakoby and S.Y. Shvartsman, Pattern formation by a moving morphogen source, \textit{Phys Biol }\textbf{8} (2011) 045003.
        \copyright IOP Publishing. Reproduced with permission. All rights reserved.}
        \end{center}
\end{figure}

\section{MATHEMATICAL FRAMEWORK}

\subsection{GENERAL SETTING}

One way to model the growth of a manifold (for instance representing the \textit{Drosophila} egg chamber) is to consider that it is provoked by a signal diffused from a source. 
We assume that the manifold is embedded in $\R^3$ and parametrized 
by the two-dimensional sphere $S^2$.
We indicate by $r$ the manifold radius, which  is a function of time $t$ 
(and of the $S^2$ coordinates) and evolves
depending on the signal $s$:
\begin{equation}\label{eq:setr}
\partial_t r = F(s).
\end{equation}
The signal $s$ also evolves with time following the reaction-diffusion equation: 
\begin{equation}\label{eq:sets}
\partial_t s = \Delta_r s + R_s(s),
\end{equation}
where $\Delta_r$ denotes the Laplace-Beltrami operator on the surface given by $r$.
Other signals may evolve on the manifold according to reaction-diffusion equations
but not inducing growth. They can be modeled by:
\begin{equation}\label{eq:setm}
\partial_t m_i = \Delta_r m_i + R_i(m_i).
\end{equation}

\subsection{TOY MODEL}

We simplify the problem described in \eqref{eq:setr}, \eqref{eq:sets} and \eqref{eq:setm} by focusing on a one-dimensional manifold and by neglecting the proteins $m_i$ that do not contribute to the growth of the cell membrane. 
Given an angle variable $\theta\in S^1$, we characterize the position
of the membrane by a function $r=r(t,\theta)$ representing the radius. 
Moreover, we consider that $F=\text{Id}$, so that the signal $s = s(t,\theta)$ directly pushes the manifold to grow in its radial direction.
The dynamics of $s$ is given by the heat equation on the manifold. We neglect the reaction term, compensating for it by allowing the signal $s$ to become negative. 
We introduce a control $u$, which sets the value of $s$ at the point $\theta = \pi$. Biologically, this corresponds to the point at which the nucleus sends the growing signal to the boundary.
Hence, the dynamics satisfies:
\begin{equation}
\begin{cases}
\partial_t r = s, \\
\partial_t s = \Delta_r s, \\
s(t, \theta=\pi) = u(t).
\end{cases}
\end{equation}

We now assume that the initial configurations of both $r$ and $s$ are symmetric with respect to $\theta$, i.e. $r(0,-\theta) = r(0,\theta)$
 and, similarly, $s(0,-\theta) = s(0,\theta)$. The simplest example is $r(0, \theta) = 1$ and $s(0,\theta) = 0$, i.e. a round cell and a zero
signal on it. One can easily prove by that, for any choice of the control $u(t)$, both $r$ and $s$ stay symmetric. Indeed, using the explicit expression of the Laplace-Beltrami operator \eqref{eq:devpde}, we prove that $(s(t,\theta), r(t, \theta))$ and $((s(t,-\theta), r(t, -\theta))$ solve the same differential system. Since the two couples have identical initial conditions, by uniqueness of solution we deduce that they are equal and thus symmetric.

Since $s$ is the solution of a heat equation, it is a $\mathcal{C}^\infty$ function far from $\theta = \pi$ for all time. As a consequence, symmetry also implies 
$\partial_{\theta} s (t,0)= 0$ for all t. Hence, we reduce our study to the half-circle 
$\theta\in [0,\pi]$ and
consider the following dynamics:
\begin{equation}
\begin{cases}
\partial_t r = s, \\
\partial_t s = \Delta_r s, \\
s(t, \theta=\pi) = u(t), \\
\partial_\theta s(t,\theta=0)=0.
\end{cases}
\end{equation}

We now study the Riemannian structure on the cell induced by a shape $r$. As already stated, $s$ is $\mathcal{C}^\infty$
except in 0, since its value there depends on $u(t)$. Assuming that the choice of $u$ implies that $s$ is $\mathcal{C}^\infty$ at 0 too,
we have that $r$ is a $\mathcal{C}^\infty$ function too. Consider the coordinate $\theta$ on the circle, and observe that a displacement
$\partial_\theta$ on the coordinate induces a displacement in the $r$ variable that can be estimated by $\sqrt{r^2+r_\theta^2} \partial_\theta$, where $r_\theta=\partial_\theta r$ is the derivative of $r$ with respect to $\theta$.
 The estimate is due to a simple geometric first-order
estimate of the length of the curve $r(\theta)$.
As a consequence, one can define the following metric on $S^1$:
\begin{equation}
g_\theta \text{ is bilinear and satisfies: } g_\theta(\partial_\theta,\partial_\theta)=r^2(\theta) + r_\theta^2(\theta).
\end{equation}

This uniquely defines the metric on $S^1$. It is also clear that the inverse of the metric satisfies $g^\theta(d\theta,d\theta)=\frac{1}{r^2(\theta) + r_\theta^2(\theta)}$.
 Such an operator is never zero since the radius is supposed to be positive for all $\theta$.
Then, a direct computation gives the explicit expression of the Laplace-Beltrami operator $\Delta_r$. We have:
\begin{equation}
\begin{split}
\Delta_r s & =\frac{1}{\sqrt{|g_\theta|}}\partial_\theta\left( \sqrt{|g_\theta|}g^\theta\partial_\theta s\right) \\
& = \frac{1}{\sqrt{r^2 + r_\theta^2}} \partial_\theta \left(\frac{1}{\sqrt{r^2 + r_\theta^2}}\partial_\theta s\right) \\
& =  \frac{1}{r^2 + r_\theta^2}\partial_\theta^2 s - \frac{r r_\theta + r_\theta \partial_\theta^2 r}{(r^2+r_\theta^2)^2} \partial_\theta s.
\end{split}
\end{equation}

Hence the system we want to study is the following:

\begin{equation}\label{eq:devpde}
\begin{cases}
\partial_t r = s, \\
\partial_t s =  \frac{1}{r^2 + r_\theta^2}\partial_\theta^2 s - \frac{r r_\theta + r_\theta \partial_\theta^2 r}{(r^2+r_\theta^2)^2} \partial_\theta s, \\
s(t, \theta=\pi) = u(t), \\
\partial_\theta s(t,\theta=0)=0.
\end{cases}
\end{equation}

We want to prove controllability for system \eqref{eq:devpde} in a specific case, that is to find a control $u$ that
drives a (symmetric) cell shape to another (symmetric) cell shape in a given time interval $[0, T]$, together
with having a signal $s$ that is zero at the initial and final times.
In mathematical terms, we consider initial and final configurations $r_0$, $r_1$ and a time $T > 0$. We want
to find a control $u : [0; T] \rightarrow \mathbb{R}$ such that the unique solution of \eqref{eq:devpde} with $r(t = 0) = r_0$ and $s(t = 0) = 0$
satisfies $r(t = T) = r_1$ and $s(t = T) = 0$. This goal is called \textit{exact controllability}. It is known that this
goal is impossible to be achieved in general, since we already know that some configurations (for instance non-smooth
final configurations) cannot be reached with a heat equation.

Hence, we instead aim to prove \textit{approximate
controllability}, defined as follows: considering initial and final configurations $r_0$, $r_1$ and a time $T > 0$, for
every $\epsilon > 0$, we want to find a control $u : [0; T] \rightarrow \mathbb{R}$ such that the unique solution of \eqref{eq:devpde} with $r(t = 0) = r0$
and $s(t = 0) = 0$ satisfies $\|r(t = T)- r_1 \|_{L^2} < \epsilon$ and $\|s(t = T)\|_{L^2} < \epsilon$.

It was shown in \cite{Laroche} that the 1-D generalized heat equation
\begin{equation}\label{eq:rouchon}
\begin{cases}
\partial_t\phi = f(\theta)\partial_{\theta}^2\phi + g(\theta)\partial_{\theta} \phi + h(\theta)\phi, \\
\phi (t, \theta=\pi)=u(t), \\
\partial_{\theta}\phi(t, \theta=0)=0
\end{cases}
\end{equation}
is approximately controllable where $f > 0$, $g$ and $h$ are analytic functions. Moreover, \cite{Laroche} proves a stronger condition:
(approximate) motion planning or (approximate) tracking, defined as follows. Given a reference trajectory, we want to find a control such that the solution of the system \eqref{eq:rouchon} stays close to the reference
trajectory for each time. In mathematical terms, one has the following result.
\begin{theorem}
Consider a time horizon $[0, T]$ and a smooth trajectory $\bar f : [0, T] \rightarrow L^2(0, \pi)$. For every $\epsilon > 0$,
there exists $u : [0, T] \rightarrow  \R$ such that the solution of \eqref{eq:rouchon} with initial data $\bar f(0)$ satisfies $\|f(t)-\bar f(t)\|_{L^2} < \epsilon$ for
all time $t \in [0, T]$.
\end{theorem}

We use this result to prove approximate controllability of \eqref{eq:devpde}. Moreover, we will show a stronger condition,
that is approximate tracking of the $r$ variable, together with the condition $\|s(t = 0)\|_{L^2} < \epsilon$ and $\|s(t =T)\|_{L^2} < \epsilon$. Since we need analytic coefficients for the second equation of \eqref{eq:devpde}, we need a reference trajectory
that is analytic for all t, i.e. $r : [0; T] \rightarrow C^w(0, \pi)$, together with smoothness with respect to $t$.
We can prove the following main theorem.
\begin{theorem}\label{th:rouchon}
Let $\bar r : [0, T] \rightarrow C^w(0,\pi)$ be a reference trajectory. Then for all $\epsilon > 0$, there exists a control
$u : [0, T] \rightarrow \R$ such that the unique solution of \eqref{eq:devpde} with $r(t = 0) = \bar r(t = 0)$ and $s(t = 0) = 0$ satisfies
$\|r(t) -\bar r(t)\|_{L^2} < \epsilon$ for all $t \in [0, T]$.
\end{theorem}

\section{EQUILIBRIA}
We look for equilibria of the form: $u(t)= u_e$ and $s(t, \theta)= s_e(\theta)$, that solves the system:
\begin{equation}\label{eq:equi}
\begin{cases}
\partial_t r_e = s_e, \\
\partial_\theta^2 s_e = \frac{r_e \partial_\theta r_e + \partial_\theta r_e \partial_\theta^2 r_e}{(r_e)^2+(\partial_\theta r_e)^2} \partial_\theta s_e, \\
s_e(\theta=\pi) = u_e, \\
\partial_\theta s_e(\theta=0)=0.
\end{cases}
\end{equation}
From \eqref{eq:devpde}, we deduce that for all $\theta$, $r_e(t,\theta)$ is a linear function of $s_e(\theta)$: 
\begin{equation}\label{eq:rlin}
r_e(t,\theta)=s_e(\theta) t+ r_0(\theta),
\end{equation}
where $r_0(\theta):=r_e(0,\theta)$.
One obvious possible equilibrium is obtained when there is no control, i.e. for a zero signal (since $s_e$ then solves a Laplace equation with the boundary condition $s_e(\theta=\pi)=0$). One gets: 
$$
\begin{cases}
u_e=0, \\
s_e(\theta) = 0 \quad \text{for all } \theta \in [0,\pi], \\
r_e(t,\theta)=r_0(\theta) \quad \text{for all } t\in [0,T], \text{ for all } \theta \in [0,\pi].
\end{cases}
$$
Hence, if $s_e$ and $u_e$ are at an equilibrium such that $u_e=0$, there is no signal and the radius is constant in time. 

On the other hand, if $u_e>0$, then $s_e$ solves a Laplace-type equation with a non-zero Dirichlet boundary condition at $\theta=\pi$, so $s_e(\theta)>0$ for all $\theta\in [0,\pi]$. Hence $r_e(t,\theta)$ grows linearly with time and does not reach an equilibrium. We instead look for an equilibrium in the shape of the membrane, by defining $\rho_e(t, \theta)=\frac{r_e(t,\theta)}{r_e(t,\theta=\pi)}$ (notice that this is possible since $r_e\neq 0$).
Then $\rho_e$ is constant in time if $\partial_t \rho_e = 0$, which gives: 
$$
\partial_t r_e(t,\theta) r_e(t,\pi) - \partial_t r_e(t,\pi) r_e(t,\theta)=0.
$$
Since $\partial_t r_e(t,\theta)=s_e(\theta)$, we get:
$$
\partial_t\rho_e(t,\theta)=0 \iff \frac{r_e(t,\theta)}{r_e(t,\pi)}=\frac{s_e(\theta)}{s_e(\pi)}.
$$
This means that at each time $t$, the membrane $r_e$ is a dilation of the signal $s_e$. 
In particular, at $t=0$, $r_0(\theta)=\frac{r_0(\pi)}{s_e(\pi)}s_e(\theta)$ for all $\theta$. Hence from \eqref{eq:rlin} we get: $r_e(t,\theta)=s_e(\theta)(t+\frac{r_0(\pi)}{s_e(\pi)})$.
Since $s_e(\theta)$ and $r_e(t,\theta)$ are proportional, the second equation of \eqref{eq:equi} becomes: 
$$
\partial_\theta^2 s_e = \frac{s_e \partial_\theta s_e + \partial_\theta s_e \partial_\theta^2 s_e}{(s_e)^2+(\partial_\theta s_e)^2} \partial_\theta s_e,
$$
which, after simplification, gives:
$$
 s_e \partial_\theta^2s_e = (\partial_\theta s_e)^2.
$$
One solution to this nonlinear differential equation is the constant signal 
$s_e(\theta)=u_e,$
where $s_e$ satisfies both the Neumann and Dirichlet boundary conditions prescribed in \eqref{eq:equi}.

We relax our conditions and look for a solution $s_e$ that satisfies $s_e(\pi)=u_e$ but not $\partial_\theta s_e(0)=0$. In particular, if we suppose that $\partial s_e (\theta) \neq 0$ for all $\theta\in [0,T]$, we can write:
$$
\frac{\partial_\theta^2 s_e}{\partial_\theta s_e}= \frac{\partial_\theta s_e}{s_e}.
$$
Then $\partial_\theta(\ln(\partial_\theta s_e))=\partial_\theta (\ln(s_e))$, so we get:
$
s_e(\theta)=u_e e^{\lambda (\theta-\pi)},
$
where $\lambda$ is a constant.
Notice that then we can bring $\partial_\theta s_e (0)=u_e \lambda e^{-\lambda \pi}$ arbitrarily close to zero by choosing $\lambda$, so we partially recover the original Neumann boundary condition.

\section{SIMULATIONS}

We simulate diffusion of the signal by discretizing the second equation of system \eqref{eq:devpde} using Finite Differences, supplemented by a Neumann boundary condition at angle $\theta=0$ ($\partial_\theta s(t,0)=0$) and a Dirichlet boundary condition at angle $\theta=\pi$ ($s(t,\pi) = u(t)$). Then the radius of the manifold at each time-step is obtained by simple integration of the signal. 

\subsection{Comparison of diffusion on a static vs growing manifold}

We run simulations for a constant control $u_1\equiv 1$, an initial signal 
$s_0(\theta)= 0$ and an initial radius $r_0(\theta)=1$ for all $\theta\in [0,\pi]$. We notice that $s$ reaches an equilibrium after time $t=2$. After that point, the radius grows in a linear way, i.e. $\rho(t)=\text{const}$. See Figure \ref{fig:uconst}.

We then turn our attention to the comparison with the case in which
we neglect the growth of the manifold (this would correspond
to an egg chamber of constant size). In this case, taking as initial condition
a circle, the radius $r$ is constant both w.r.t. time and the $\theta$ variable,
thus $r\equiv 1$ and $r_{\theta}\equiv 0$. Plugging this information
into equation (\ref{eq:devpde}), the Laplace-Beltrami operator
reduces to standard diffusion and we get the following system: 
\begin{equation}\label{eq:standard}
\begin{cases}
\partial_t r = s, \\
\partial_t s =  \partial_\theta^2 s \\
s(t, \theta=\pi) = u(t), \\
\partial_\theta s(t,\theta=0)=0.
\end{cases}
\end{equation}
The simulations for a constant control $u\equiv 1$ are very different from those obtained by using the system (\ref{eq:devpde}): Figure \ref{fig:stand} 
shows the evolution of the signal and the radius with constant control for system.
The signal $s$ reaches an equilibrium $s(t,\theta)=1$, which means that the growth of the radius tends to be uniform with respect to the angle $\theta$. 
Therefore, as expected, neglecting the growth of the manifold
generates uniform growth and, in the biological system, would give rise
to spherical egg chambers opposed to the spheroidal ones observed in nature. 

\begin{figure}[h!]
        \begin{center}
                \includegraphics[trim=3cm 0cm 3cm 0cm, clip=true, scale=0.2]                                {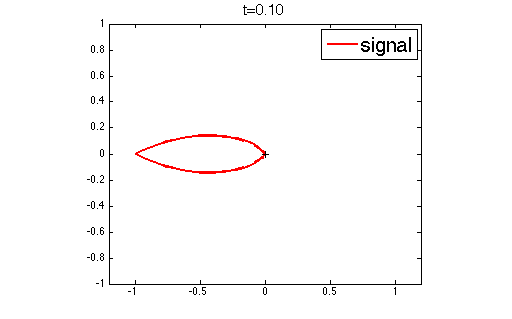} 
                \hspace{0cm}
                \includegraphics[trim=2cm 0cm 2cm 0.5cm, clip=true, scale=0.2]{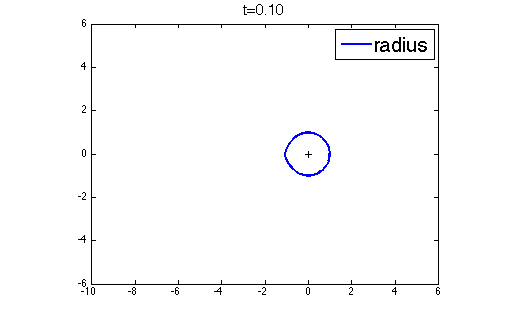}             
\\
                 \includegraphics[trim=3cm 0cm 3cm 0cm, clip=true, scale=0.2]                                {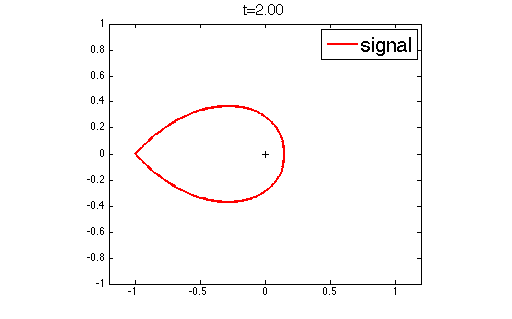} 
                \hspace{0cm}
                \includegraphics[trim=2cm 0cm 2cm 0.5cm, clip=true, scale=0.2]{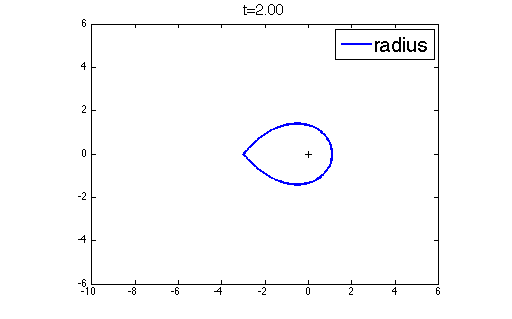} 
\\                
                \includegraphics[trim=3cm 0cm 3cm 0cm, clip=true, scale=0.2]                                 {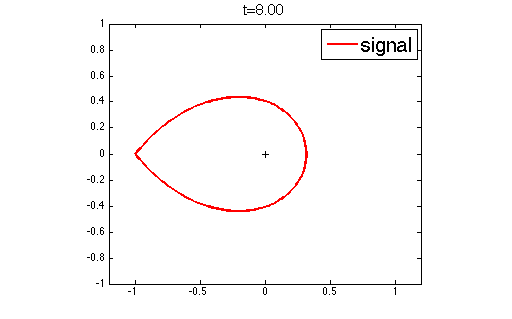} 
                \hspace{0cm}
                \includegraphics[trim=2cm 0cm 2cm 0.5cm, clip=true, scale=0.2]{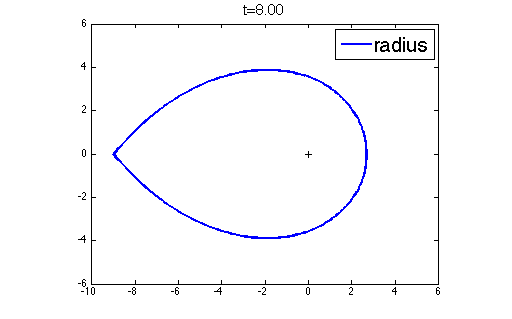} 
                \label{fig:uconst}

        \caption{Signal $s$ (left) and radius $r$ (right) for a constant control $u\equiv 1$ at times $t=0.1$, $t=2$ and $t=8$. The source correponds to the
        angle $\theta=\pi$, so in the signal picture it is located on the left end of the equator line corresponding to coordinates $(-1,0)$.}\label{fig:uconst}
        \end{center}
\end{figure}

\begin{figure}[h!]
        \begin{center}
                \includegraphics[trim=3cm 0cm 3cm 0cm, clip=true, scale=0.2]{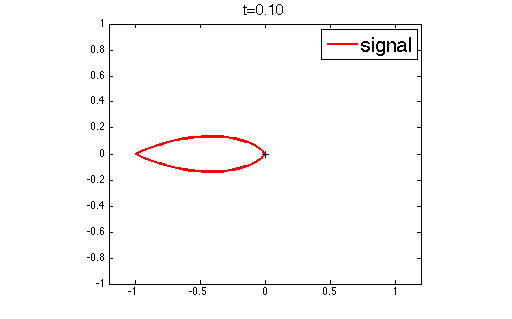} 
                \hspace{0cm}
                \includegraphics[trim=2cm 0cm 0cm 0cm, clip=true, scale=0.2]{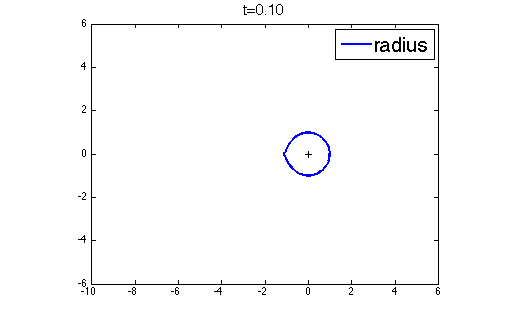}             
\\
                \includegraphics[trim=3cm 0cm 3cm 0cm, clip=true, scale=0.2]{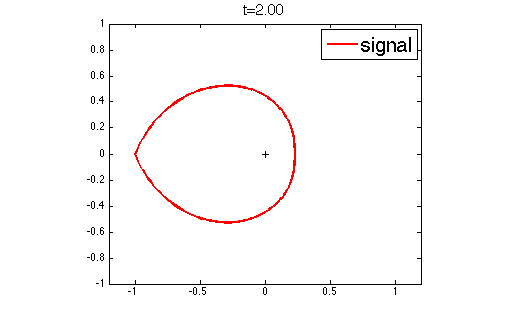} 
                \hspace{0cm}
                \includegraphics[trim=2cm 0cm 0cm 0cm, clip=true, scale=0.2]{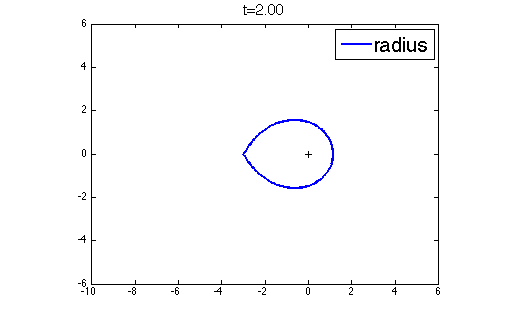} 
                \\
                 \includegraphics[trim=3cm 0cm 3cm 0cm, clip=true, scale=0.2]{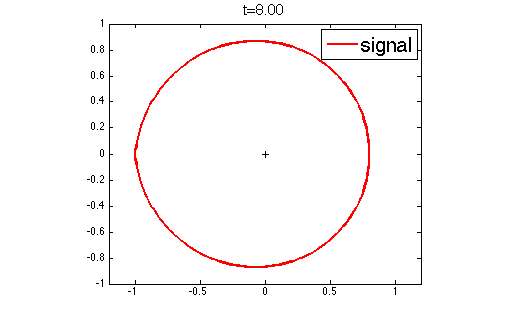} 
                \hspace{0cm}
                \includegraphics[trim=2cm 0cm 0cm 0cm, clip=true, scale=0.2]{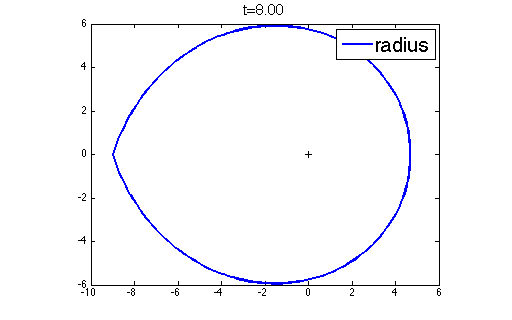} 

        \caption{Signal $s$ (left) and radius $r$ (right) for a constant control $u\equiv 1$ at times $t=1$, $t=2$, $t=5$ and $t=8$.}\label{fig:stand}
        \end{center}
\end{figure}

\subsection{Single source}
\begin{figure}[h!]
        \begin{center}
                \includegraphics[trim=0cm 0cm 0cm 0cm, clip=true, scale=0.3]{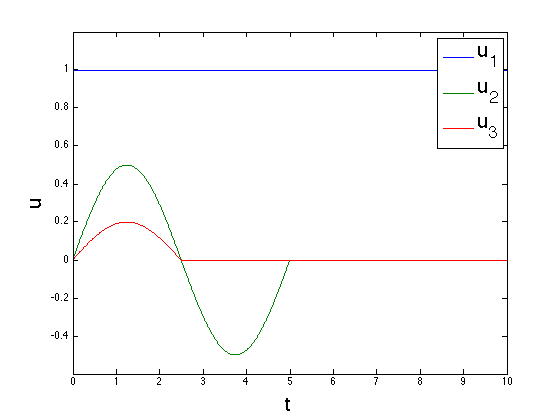} 
        \caption{Control functions $u_1$, $u_2$ and $u_3$}\label{fig:controls}
        \end{center}
\end{figure}

\begin{figure}[h!]
        \begin{center}
                \includegraphics[trim=2cm 0cm 2cm 0cm, clip=true, scale=0.25]{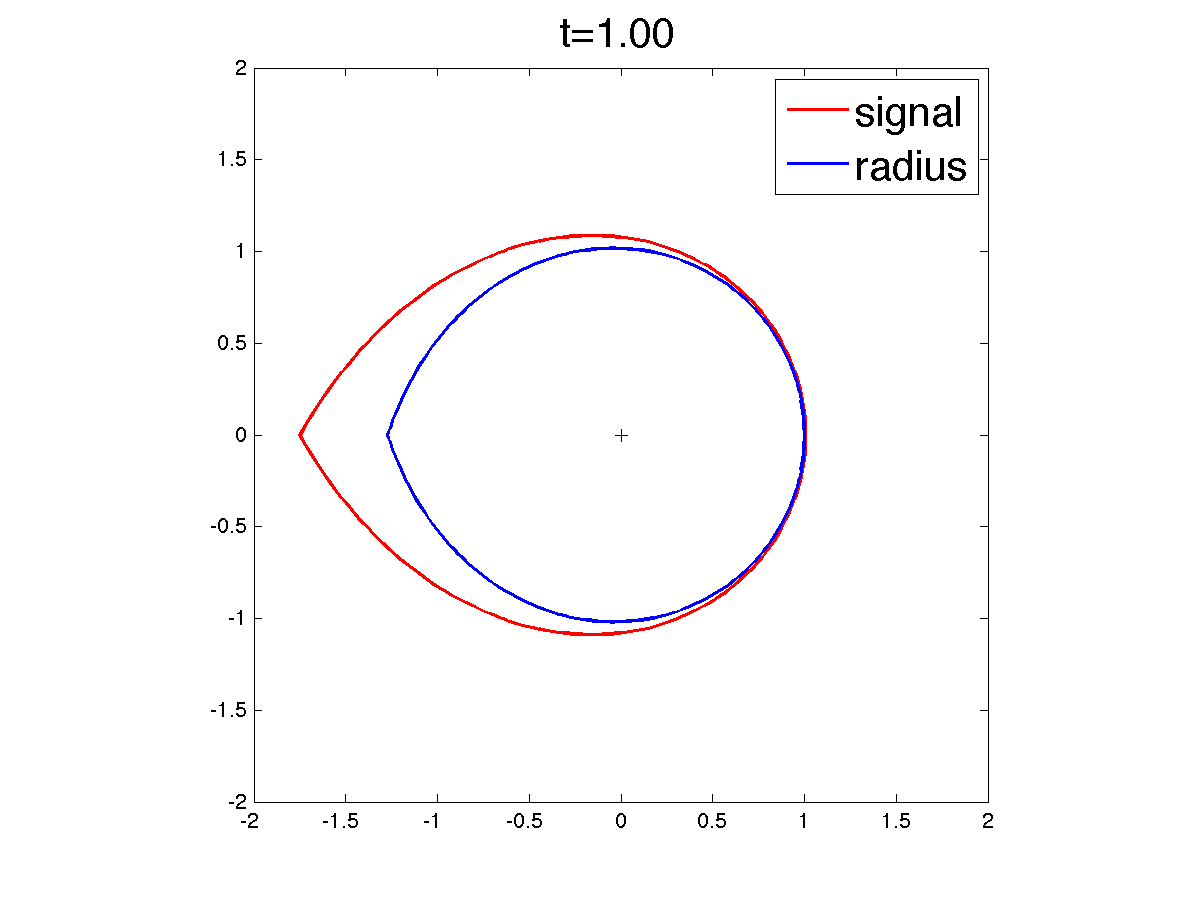} 
                \includegraphics[trim=2cm 0cm 2cm 0cm, clip=true, scale=0.25]{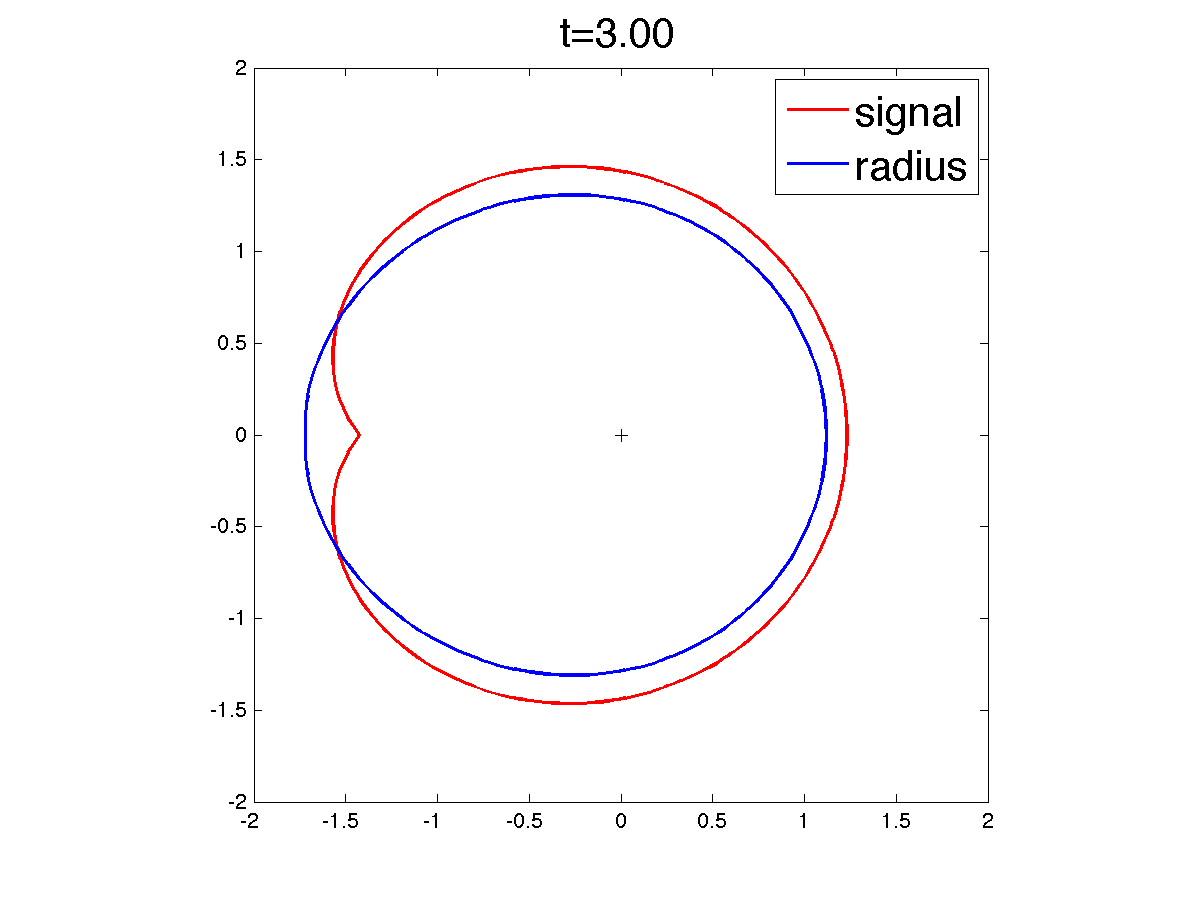}             

                \includegraphics[trim=2cm 0cm 2cm 0cm, clip=true, scale=0.25]{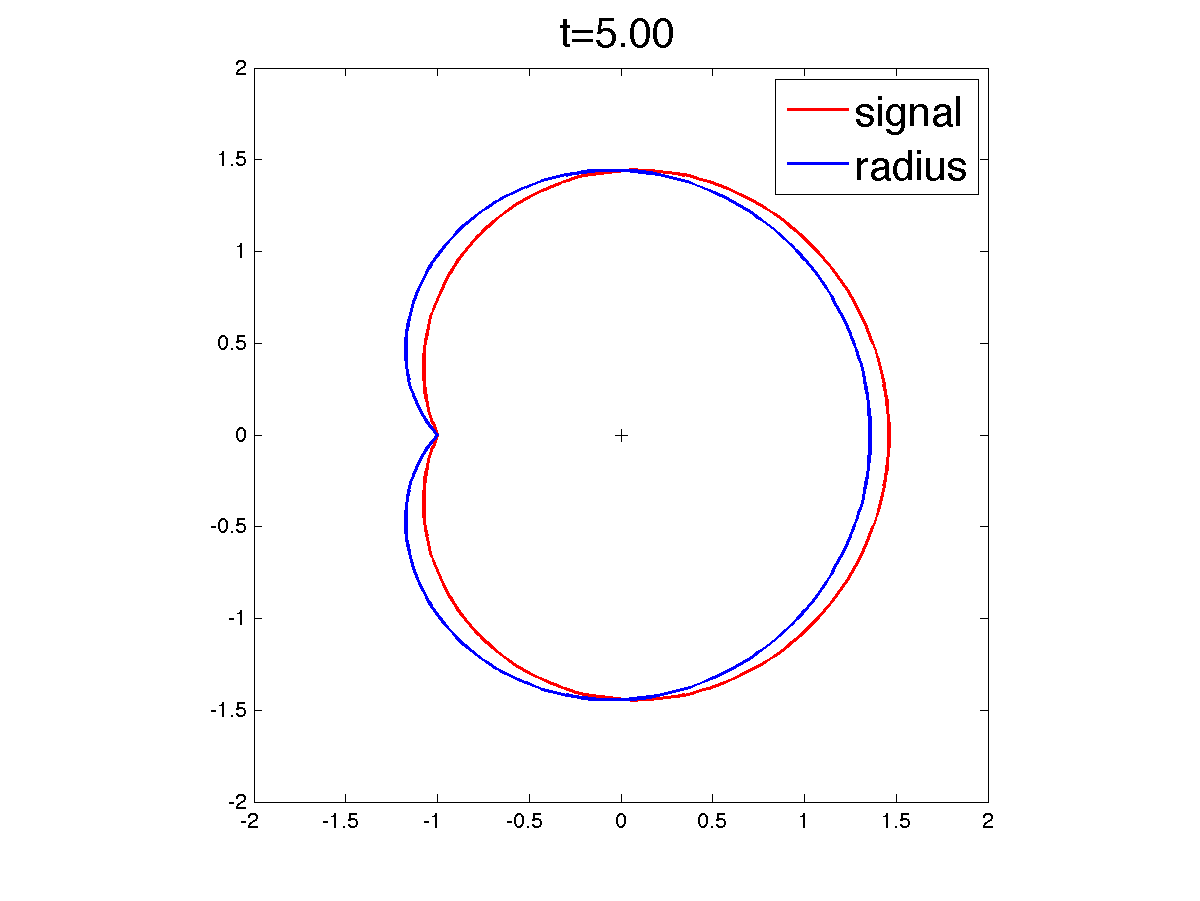} 
                \includegraphics[trim=2cm 0cm 2cm 0cm, clip=true, scale=0.25]{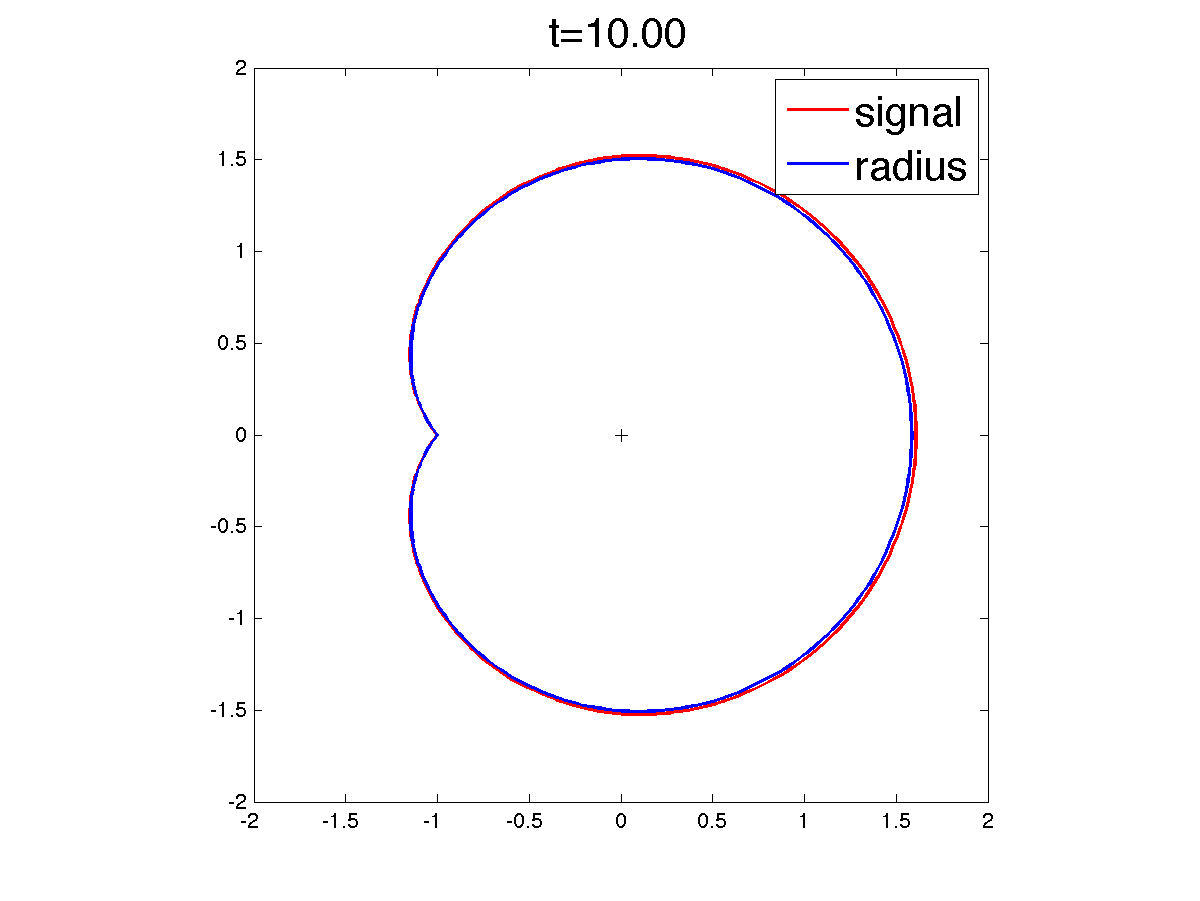} 
        \caption{Radius $r$ (in blue) and signal $s$ (plotted as $r+s$ in red) for a control $u=u_2$ at times $t=1, 3, 5$, and $10$.}\label{fig:rightgrowth}
        \end{center}
\end{figure}

A source placed on the first axis (at angle $\theta = \pi$) allows us to control the diameter of the manifold along the same axis. In Figure \ref{fig:uconst}, the manifold is stretched along the first axis direction at final time, with an emphasis on the left side, i.e. $r(T,\pi)>r(T,0)$.
Using the source to impose negative values of the signal (which has a mathematical meaning but not a biological one), we can control the final shape of the manifold to achieve $r(T,\pi)<r(T,0)$. In order to do that we set the control as:
\[
u(t)=\left\{
\begin{array}{ll}
0.5\cdot sin(\omega t) & t\in [0,5]\\
0 & t\in ]5,10]\\
\end{array}
\right. ,
\]
where $\omega=\frac{2\pi}{5}$ so that we obtain a complete sinusoidal oscillation
up to time $5$ then the signal is vanishing (which coincides with control $u_2$
depicted in Figure \ref{fig:controls}). 
The final result is a apple shape manifold
with pitch located at the signal source point, see Figure \ref{fig:rightgrowth}.
To better visualize the relationship between the signal and the shape
we visualized the signal on the manifold itself, so for positive values the
signal will be outside the manifold and inside for negative ones. 

\begin{figure}[h!]
        \begin{center}
                \includegraphics[trim=2cm 0cm 2cm 0cm, clip=true, scale=0.25]{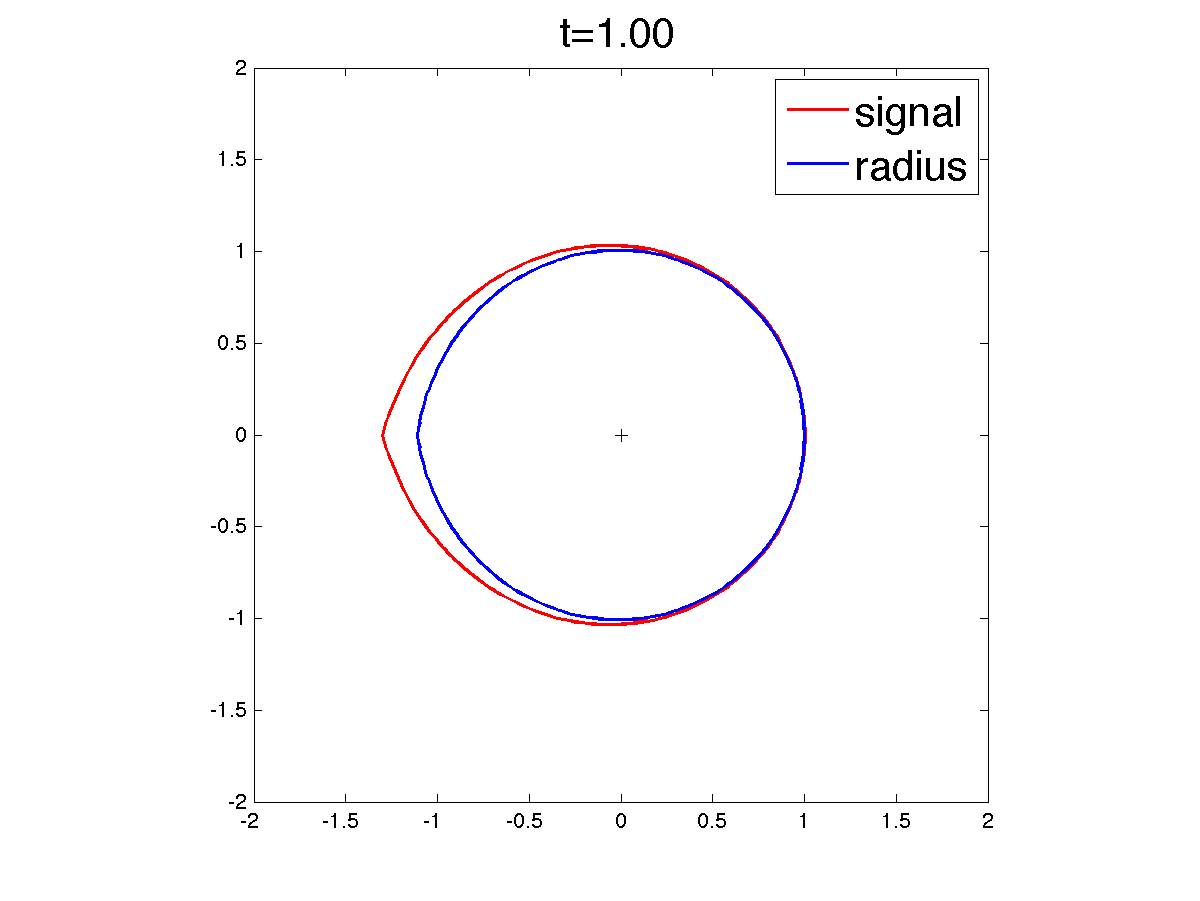} 
                \includegraphics[trim=2cm 0cm 2cm 0cm, clip=true, scale=0.25]{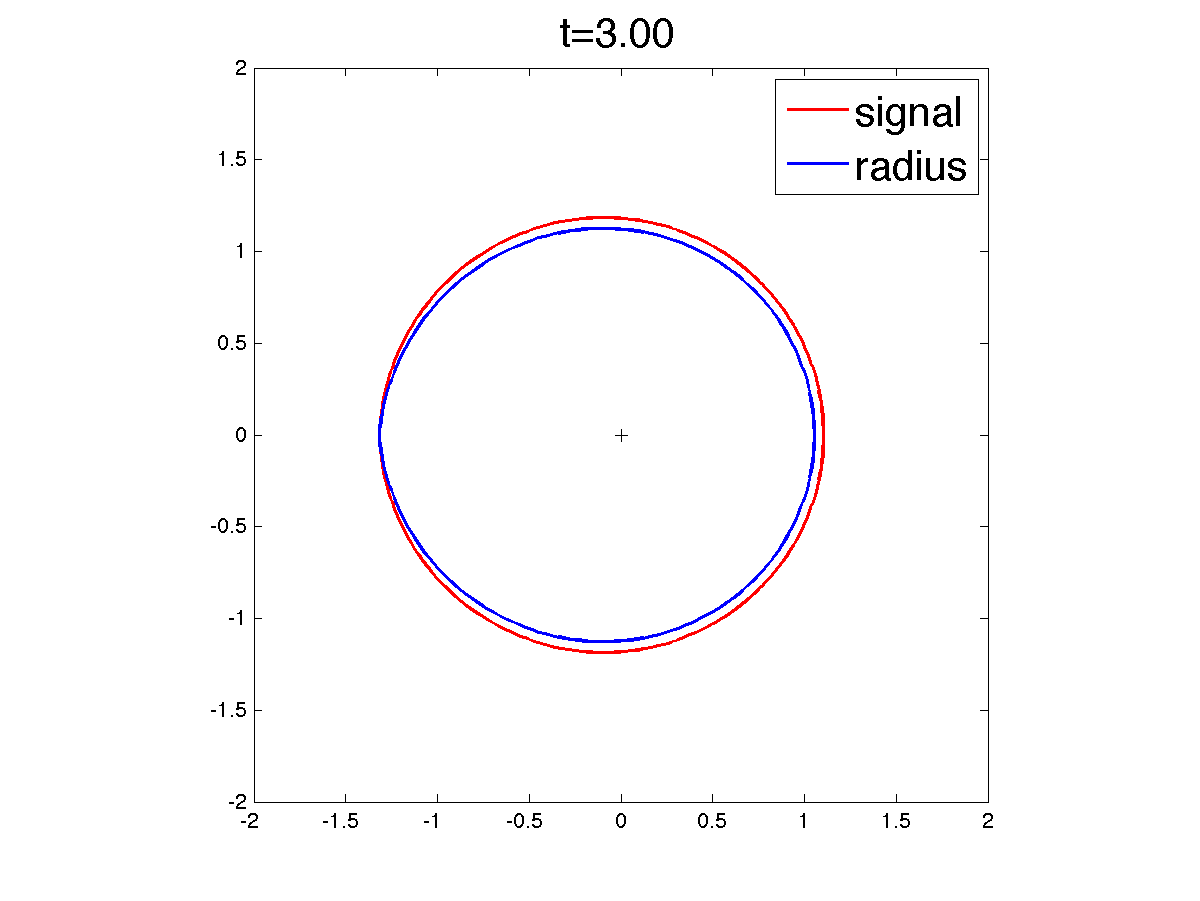}             

                \includegraphics[trim=2cm 0cm 2cm 0cm, clip=true, scale=0.25]{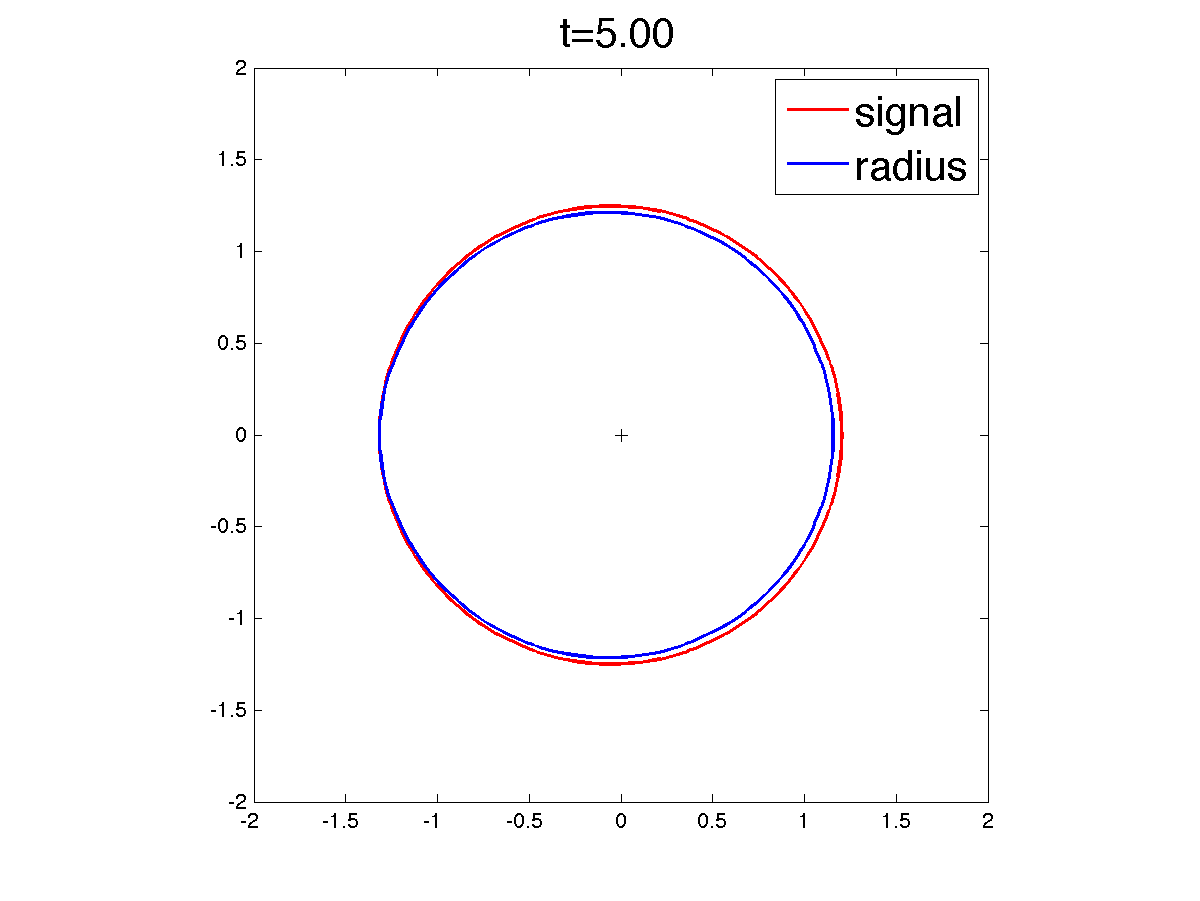} 
                \includegraphics[trim=2cm 0cm 2cm 0cm, clip=true, scale=0.25]{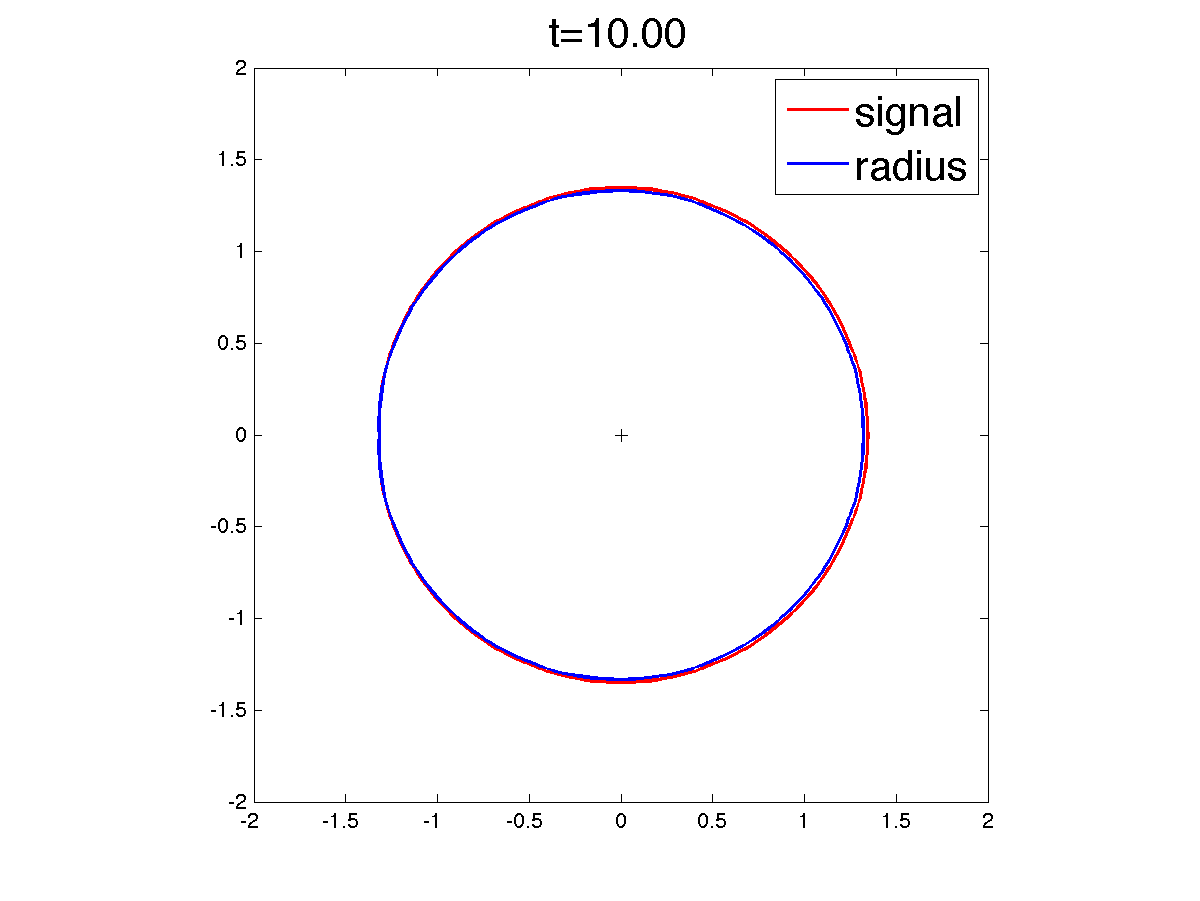} 
        \caption{Radius $r$ (in blue) and signal $s$ (plotted as $r+s$ in red) for a control $u=u_3$ at times $t=1, 3,  5$ and $10$.}\label{fig:circle}
        \end{center}
\end{figure}

Using a single source it is also possible to induce an homogeneous growth
along all directions, but with time-dependent signals. We first give an impulse
and then turn off the signal. Define the control by:
\begin{equation}\label{eq:impulse}
u(t)=\left\{
\begin{array}{ll}
0.2\cdot sin(\omega t) & t\in [0,2.5]\\
0 & t\in ]5,10]\\
\end{array}
\right. ,
\end{equation}
where $\omega=\frac{2\pi}{5}$, so that the half sinusoidal oscillation gives
an always positive signal (this correspond also to the control $u_3$
depicted in Figure \ref{fig:controls}). 
The final shape is close to that of a circle, but with a larger radius than that at initial time (see Figure \ref{fig:circle}).

\subsection{Double source}
\begin{figure}[h!]
        \begin{center}
                \includegraphics[trim=2cm 0cm 2cm 0cm, clip=true, scale=0.25]{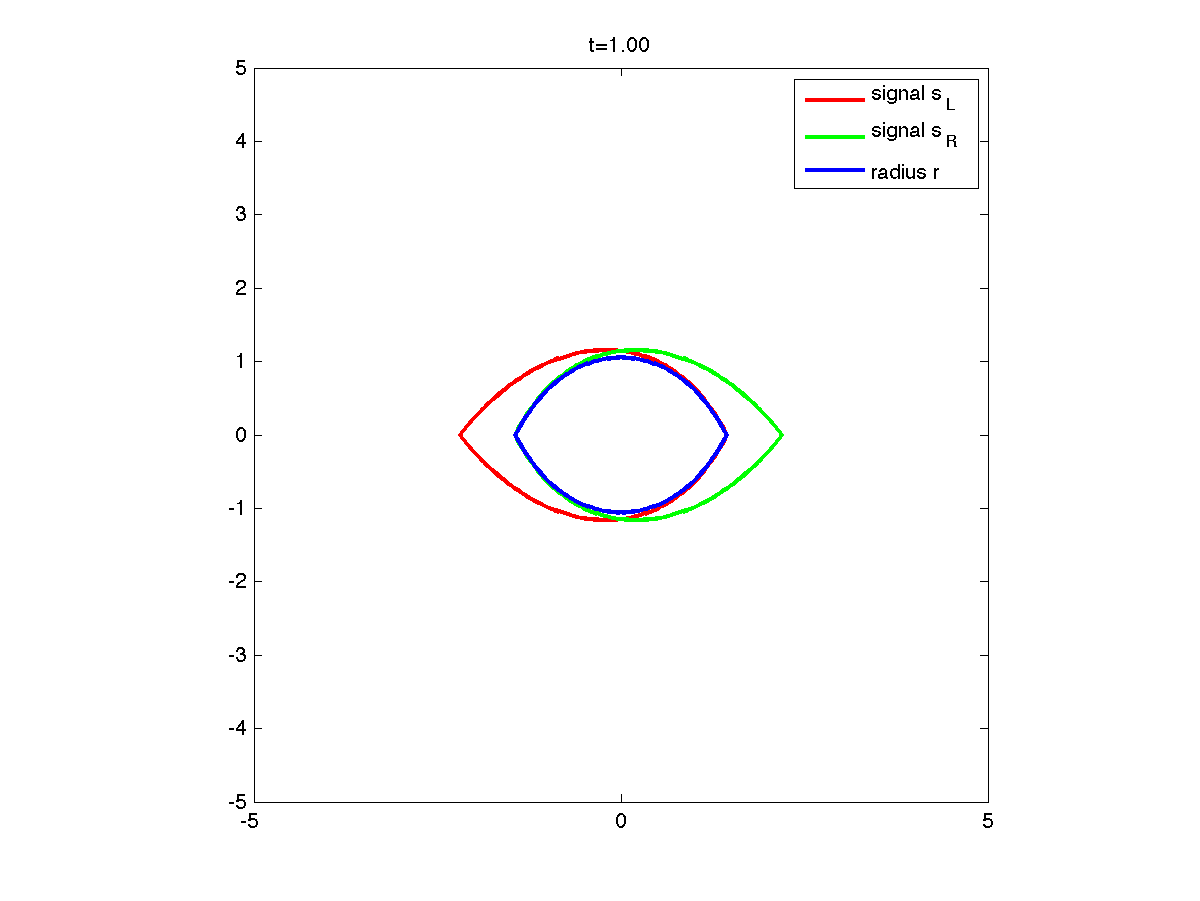} 
                \includegraphics[trim=2cm 0cm 2cm 0cm, clip=true, scale=0.25]{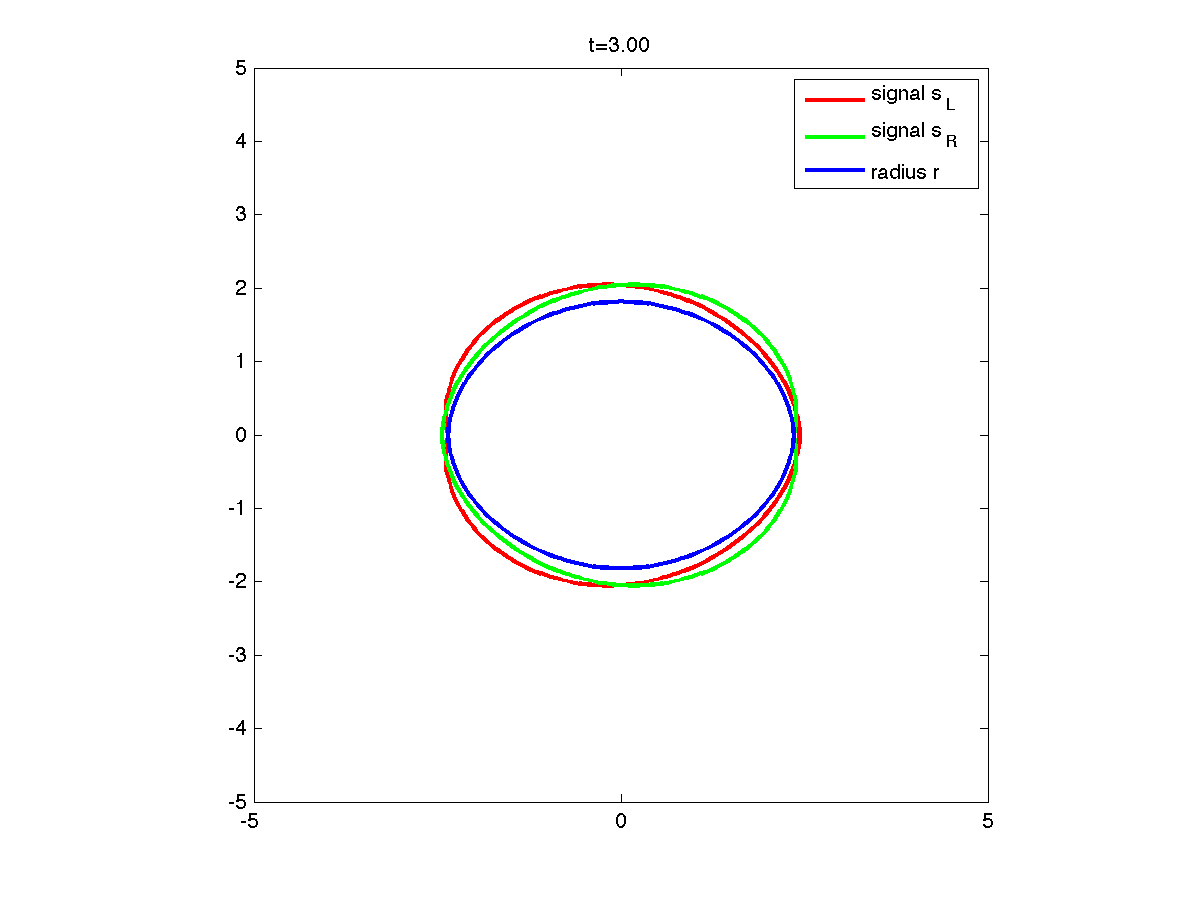}             

                \includegraphics[trim=2cm 0cm 2cm 0cm, clip=true, scale=0.25]{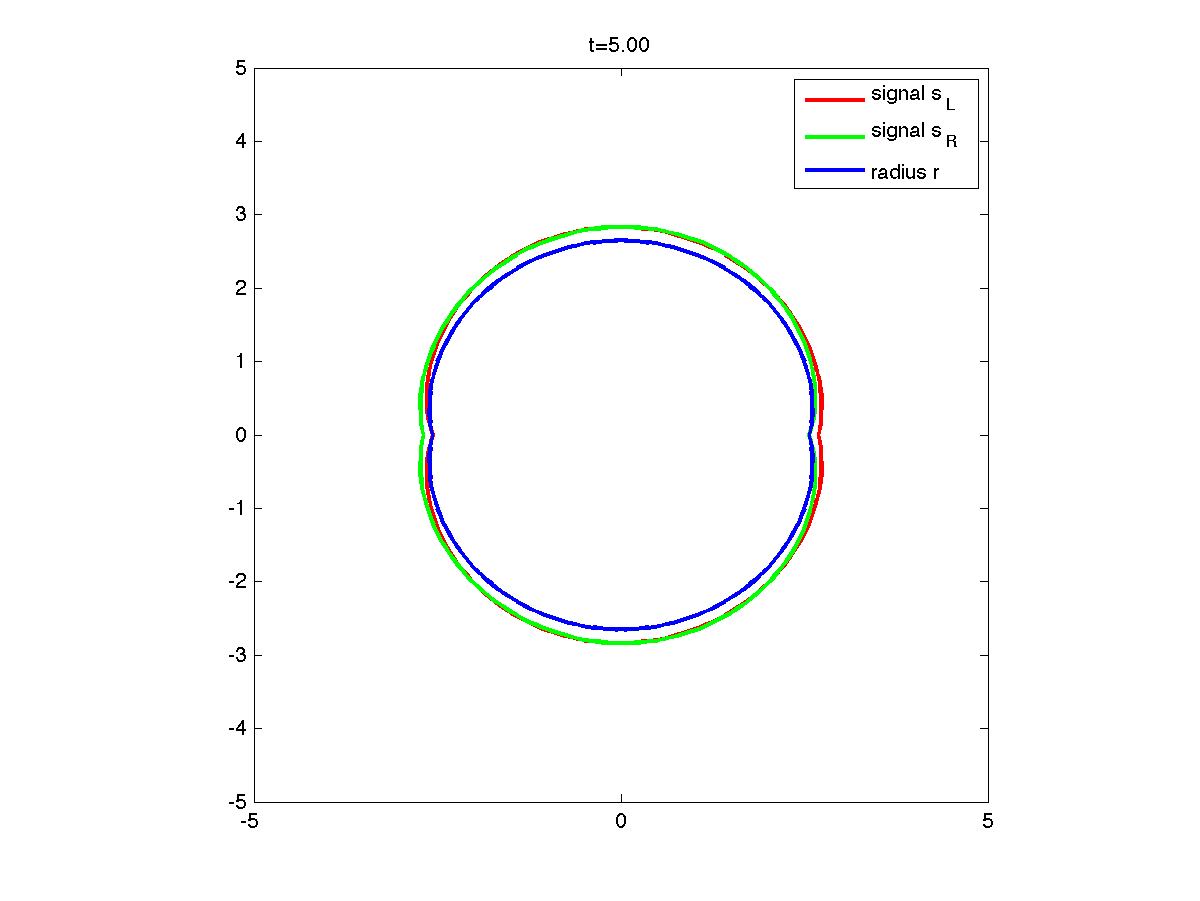} 
                \includegraphics[trim=2cm 0cm 2cm 0cm, clip=true, scale=0.25]{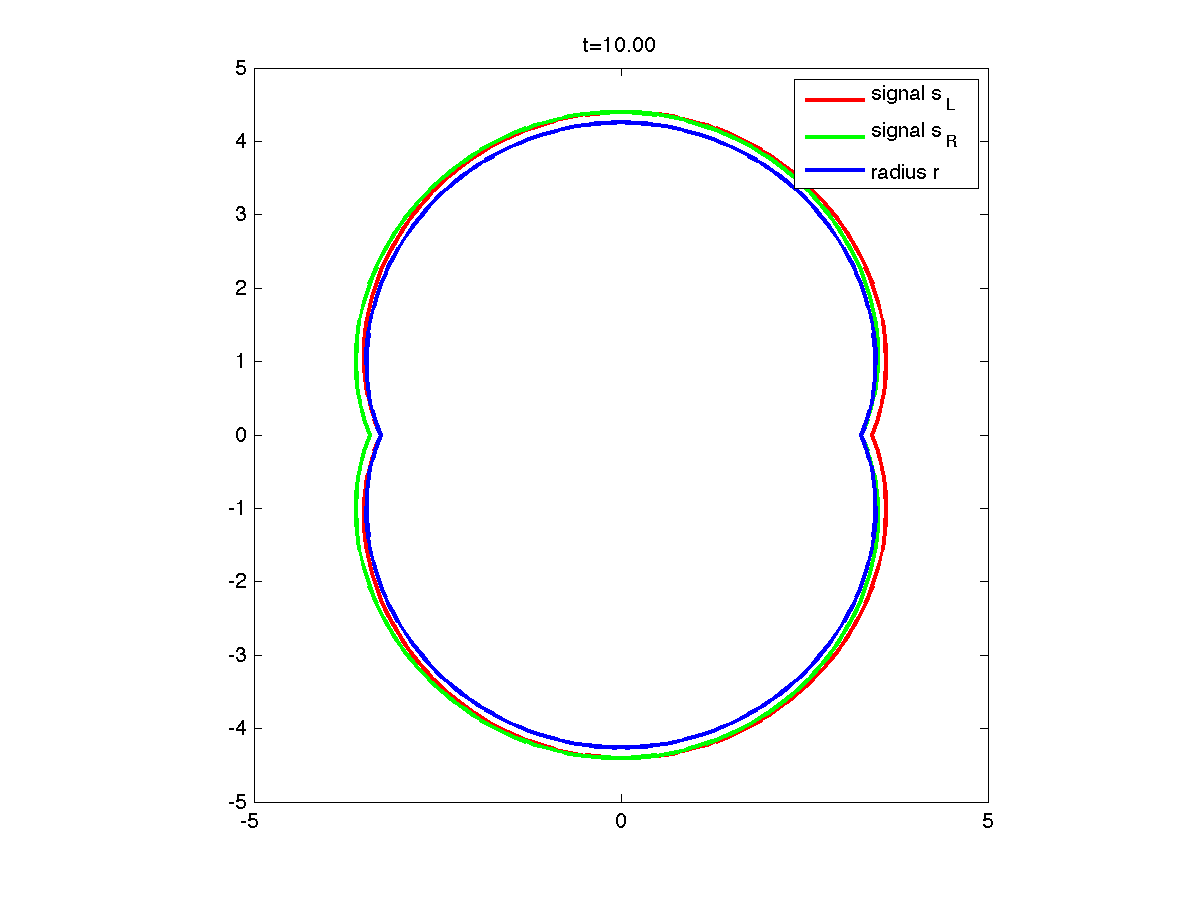} 
        \caption{Radius $r$ (in blue) and signals $s_L$ (plotted as $r+s_L$ in red) and $s_R$ (plotted as $r+s_R$ in green) for controls $u_L=u_R=u_3$ at times $t=1$, $t=3$, $t=5$ and $t=10$.}\label{fig:2sources}
        \end{center}
\end{figure}

As observed above, a single static source allows us to control the radii $r(T,0)$ and $r(T,\pi)$, i.e. the horizontal growth. In order to achieve a larger growth
along the vertical axis, we consider a system with double source: one
locate at angle $\theta=0$ and the second (as before) at angle $\theta=\pi$.
We obtain the system:
\begin{equation}\label{eq:2sources}
\begin{cases}
\partial_t r =s_L+s_R, \\
\partial_t s_L =  \frac{1}{r^2 + r_\theta^2}\partial_\theta^2 s_L - \frac{r r_\theta + r_\theta \partial_\theta^2 r}{(r^2+r_\theta^2)^2} \partial_\theta s_L \\
\partial_t s_R =  \frac{1}{r^2 + r_\theta^2}\partial_\theta^2 s_R - \frac{r r_\theta + r_\theta \partial_\theta^2 r}{(r^2+r_\theta^2)^2} \partial_\theta s_R \\
s_L(t, \theta=\pi) = u_L(t), \quad  s_R(t, \theta=0) = u_R(t),\\
\partial_\theta s_L(t,\theta=0)=0, \quad \partial_\theta s_R(t,\theta=\pi)=0.
\end{cases}
\end{equation}
If we use the control given by formula (\ref{eq:impulse}) for both sources,
we obtain a final manifold stretched more in the vertical direction, i.e. $r(T,\pi/2)>r(T,0)=r(T,\pi)$, see Figure \ref{fig:2sources}.

\section*{CONCLUSION} 
In this paper we introduced a new mathematical framework, called
Developmental Partial Differential Equation (DPDE), to model
the growth of organisms induced by signaling pathways.
A DPDE consists of a couple: a time-varying manifold and a signal
evolving on the manifold. 
Inspired by the specific application to \textit{Drosophila} egg chamber development,
we consider a completely  coupled evolution where the manifold's growth is regulated by the signal and the signal diffusion 
by an operator (Laplace-Beltrami) depending on the manifold geometry.\\
We provide controllability results using flatness of the heat equation
and show simulations of resulting manifold shapes.\\
Future work will include: moving sources, general
differential operators, higher dimensions manifolds
and explicit expressions of controls for motion planning.

\addtolength{\textheight}{-11cm}   

\section*{ACKNOWLEDGMENT}

The authors acknowledge the support of the NSF Project "KI-Net", DMS Grant \# 1107444.


\end{document}